\documentclass{amsart}
\usepackage{graphicx}
\usepackage{mathpazo}

\newtheorem{thm}{Theorem}[section]

\newtheorem{lem}[thm]{Lemma}

\theoremstyle{definition}

\theoremstyle{remark}

\numberwithin{equation}{section}
\newcommand{\norm}[1]{\left\Vert#1\right\Vert}
\newcommand{\abs}[1]{\left\vert#1\right\vert}

\newcommand{\eps}{\varepsilon}

\newcommand{\ye}{y_\varepsilon}
\newcommand{\lme}{\lambda^\varepsilon}
\newcommand{\eto}{\varepsilon\to0}
\begin{document}

\title[Spectral properties of the fourth order differential operator]
{On spectral properties of the fourth order differential operator with singular coefficients}%
\author{Stepan Man'ko}%
\address{Department of Mechanics and Mathematics, Ivan Franko National University of Lviv, 1 Universytetska str., 79000 Lviv, Ukraine}%
\email{s\_\,manko@franko.lviv.ua}%

% \thanks{}%
% \subjclass{}%
% \keywords{}%

%\date{}%
%\dedicatory{}%
%\commby{}%
% ----------------------------------------------------------------
\begin{abstract}
A formal fourth order differential operator with a singular coefficient that is a linear combination of the Dirac delta-function and its derivatives is considered.
The asymptotic behavior of spectra and eigenfunctions of a family of differential operators with smooth coefficients approximating the singular coefficients is studied.
We explore how behavior of eigenvalues and eigenfunctions is influenced by singular coefficients. The limit operator is constructed and is shown to depend on a type of approximation of singular coefficients.
\end{abstract}
\maketitle
% ----------------------------------------------------------------
\section{Introduction}
Differential operators with singular coefficients appear in atomic physics, acoustics, quan\-tum mechanics, solid state physics, aerodynamics, fluid mechanics, aeroacoustics \cite{Alb2ed}, \cite{AlbKur}, \cite{DemkovOstrovskii}, \cite{GolshOri}, \cite{Yablonski},  \cite{Zavalishchin}.
An important task of the theory of differential equations is to find the minimal smoothness of coefficients, under which the equation admits a solution.
Although there are some models that are closely related to the differential operators with distributions in coefficients it is impossible to construct the theory of linear differential equations with distributional coefficients, since  the space of distributions $\mathcal{D}'(\mathbb{R}^n)$ is not an algebra with respect to the ``pointwise'' multiplication.
This raises the basic question how to interpret the differential operators with distributions in coefficients.

A lot of models are expected to be ``selfadjoint'' in the sense that  appropriate operators, describing these models must be selfadjoint in some Hilbert spaces.
Let a differential expression  $S$ correspond to such a model and let it contain distributions, supported by $x=0$, in coefficients.
In order to interpret the operator $S$ we first construct a symmetric operator $S_0$ by restricting $S$ to the set of functions vanishing at the origin along with their derivatives.
Then we consider all selfadjoint extensions of $S_0$ and choose one of these extensions as a definition of $S$.
This method goes back to the work of F. Berezin and L. Faddeev \cite{BerezinFadeev}. In some instances the set of all selfadjoint extensions of a symmetric operators
is multiparametric.
Therefore the harder question comes: how to choose
an extension that is best suited to our physical model. For some models
the proper operator can not be chosen within the selfadjoint extensions
theory, because the models contain hidden parameters. After replacing the
singular coefficient with a sequence of short-range smooth coefficients, the operator obtained in the zero-range limit, as often happens, can depend
on the type of regularization, i.e., the operator is governed by the
shape of squeezed coefficients. This shape is a hidden parameter and plays
a crucial role in the choice of a selfadjoint extension corresponding to the
physical model under consideration.

In \cite{GolovManko0}, \cite{GolovManko1} the problem how to define the one-dimensional Schr\"{o}dinger operator with the $\delta'$-potential, where $\delta'$ is the first derivative of the Dirac delta-function, was considered .
A natural approach to defining such a Hamiltonian is to approximate $\delta'$ in $\mathcal{D}'(\mathbb{R}^n)$-topology by regular
potentials and then to investigate the corresponding family of regular Schr\"{o}dinger operators.
Therefore the authors considered the family of Schr\"{o}dinger operators on
the line of the form
$$
\mathcal{H}_\eps(\alpha,\Psi)=-\frac{d^2}{dx^2}+U(x)+
\frac{\alpha}{\eps^2}\Psi\Big(\frac{x}{\eps}\Big),
$$
approaching a formal Hamiltonian $H_\alpha=-\frac{d^2}{dx^2}+U(x)+
\alpha\delta'(x)$.
Here $\eps$ is a small positive parameter, $\Psi\in C_0^\infty(-1,1)$,
$U$ is a real valued potential going to $+\infty$ as $|x|\to\infty$, and
$\alpha$ is a real coupling constant.
The map assigning a limit operator $\mathcal{H}(\alpha,\Psi)$
to each pair $(\alpha,\Psi)$ was constructed there.
The choice of $\mathcal{H}(\alpha,\Psi)$ is determined by proximity of its energy levels and pure states to those for the Hamiltonian
with regularized potentials for small $\eps$.
For almost all $\alpha$ the limit operator is just the direct sum of the Schr\"{o}dinger operators with the potential $U$ on
half-axes subject to the Dirichlet boundary condition at the origin (the \textit{nonresonant case}).
But for $\alpha$ belonging to the discrete {\it resonant set} $\Sigma_\Psi$, which is the spectrum of the
Sturm-Liouville problem $-w''+\alpha\Psi w=0$ on the interval $(-1,1)$ subject to the boundary conditions $w'(-1)=w'(1)=0$,
the operator $\mathcal{H}(\alpha,\Psi)$ acts via $\mathcal{H}(\alpha,\Psi)f=-f''+Uf$ on an appropriate set of functions satisfying the matching conditions $f(+0)=\theta_\Psi(\alpha)f(-0)$
and $\theta_\Psi(\alpha)f'(+0)=f'(-0)$ (the {\it resonant case}).
Here $\theta_\Psi(\alpha)=w_\alpha(1)/w_\alpha(-1)$, where $w_\alpha$ is an eigenfunction corresponding to the eigenvalue $\alpha\in\Sigma_\Psi$.

In \cite{Manko} the results of \cite{GolovManko0}, \cite{GolovManko1} were extended to the case of the fourth order ordinary differential operators.
An attempt was made to define the formal differential operator $A_\alpha=\frac{d^4}{dx^4}+U(x)+\alpha\delta'''(x)$.
To approximate $A_\alpha$, the family $\mathcal{A}_\eps(\alpha,\Psi)=\frac{d^4}{dx^4}+U(x)+
\frac{\alpha}{\eps^4}\Psi\big(\frac{x}{\eps}\big)$ with domain
\begin{align*}
\mathcal{D}(\mathcal{A}_\eps(\alpha,\Psi))=\{f\in W_2^4(a,b)\colon\;\; f(a)=f'(a)=0,\quad
f(b)=f'(b)=0\}
\end{align*}
was considered. Here $(a,b)$ is an interval of $\mathbb{R}$  containing the origin,  $U$ is a smooth real valued function on $[a,b]$ and $\Psi\in C_0^\infty(-1,1)$.
Asymptotic expansions for eigenvalues and eigenfunctions of $\mathcal{A}_\eps(\alpha,\Psi)$
were constructed, and therefore the limit operator $\mathcal{A}(\alpha,\Psi)$ was obtained.
Upon constructing asymptotics two different cases are distinguished: the {\it resonant case} and the {\it nonresonant} one.
In the resonant case $\alpha$ belongs to the discrete {\it resonant set} $\Sigma_\Psi\subset\mathbb{R}$, which is the spectrum of the
eigenvalue problem
\begin{equation}\label{Resonant Set}
\begin{cases}
w^{(4)}+\alpha\Psi w=0,\quad \xi\in(-1,1),\\
w''(-1)=w'''(-1)=0,\quad w''(1)=w'''(1)=0.
\end{cases}
\end{equation}
The limit operator was obtained under the additional assumptions
\begin{equation}\label{MainAssumption}
\alpha\;{\it is\; a\; simple\; eigenvalue\; of\; the\; problem}\;\eqref{Resonant Set}\; {\it and}\;
w_\alpha'(-1)w_\alpha'(1)\neq0,
\end{equation}
where $w_\alpha$ is an eigenfunction corresponding to the eigenvalue $\alpha$.
In this case $\mathcal{A}(\alpha,\Psi)$ acts via $\mathcal{A}(\alpha,\Psi)f=f^{(4)}+Uf$ on an appropriate set of functions obeying the interface conditions
$f(0)=0$, $f'(+0)=\theta_\Psi(\alpha)f'(-0)$ and $\theta_\Psi(\alpha)f''(+0)=f''(-0)$, where $\theta_\Psi(\alpha)=w'_\alpha(1)/w'_\alpha(-1)$.
In the nonresonant case, when $\alpha\notin\Sigma_\Psi$, the limit operator is the direct sum of the Dirichlet operators on $(a,0)$ and $(0,b)$ respectively.

We extend the results of \cite{Manko} to more a general perturbation of the fourth order differential operator, namely, we consider a formal differential expression
$$
\frac{d^4}{dx^4}+U(x)+\alpha \delta'''(x)+\beta \delta''(x)+
\gamma_1\delta'(x)+\gamma_2\delta(x).
$$
The investigation of the papers \cite{GolovManko0}--\cite{Manko} is based on the asymptotic analysis.
We will use techniques of \cite{GolovManko0}--\cite{Manko} to obtain an appropriate limit operator.

\subsection{Problem statement and main results}
Let $L$ stand for the differential expression $\frac{d^4}{dx^4}+U(x)$.
As before $U$ is a smooth real valued function on the
interval $[a,b]\subset\mathbb{R}$, containing the origin.
Denote by $\Psi_\eps$ the function
\begin{equation*}
    \Psi_\eps(x)=\alpha\eps^{-4}\Psi(\eps^{-1}x)
    +\beta\eps^{-3}\Phi(\eps^{-1}x)+\gamma_1\eps^{-2}\Upsilon_1(\eps^{-1}x)+
    \gamma_2\eps^{-1}\Upsilon_2(\eps^{-1}x).
\end{equation*}
Here $\Psi,\Phi,\Upsilon_1,\Upsilon_2\in C_0^\infty(-1,1)$,
$\mathrm{supp}\,\Psi=[-1,1]$, and $\alpha,\beta,\gamma_1,\gamma_2\in\mathbb{R}$ are arbitrary con\-stants.
Let us consider the eigenvalue problem
\begin{equation}\label{MainProblem}
L\ye+\Psi_\eps(x)\ye=\lme \ye,\quad x\in(a,b),\qquad
        \ye(a)=\ye'(a)=\ye(b)=\ye'(b)=0.
\end{equation}
Note that the further analysis of the problem \eqref{MainProblem} does not depend on the type of boundary conditions. Hence Dirichlet boundary conditions may be replaced
by one of the possible combinations at the endpoints $x=a$ and $x=b$ of the following boundary conditions
\begin{align*}
    y(x_0)=y'(x_0)=0,\qquad
    y(x_0)=y''(x_0)=0,\qquad
    y''(x_0)=y'''(x_0)=0.
\end{align*}
We associate with the problem \eqref{MainProblem} an operator
\begin{align*}
&\mathcal{S}_\eps(\alpha,\beta,\gamma_1,\gamma_2;\Psi,\Phi,\Upsilon_1,\Upsilon_2)=\frac{d^4}{dx^4}+U(x)+\Psi_\eps(x),\\
&    \mathcal{D}(\mathcal{S}_\eps(\alpha,\beta,\gamma_1,\gamma_2;\Psi,\Phi,\Upsilon_1,\Upsilon_2))=
    \{f\in W_2^4(a,b)\colon\\
    &\qquad \qquad \qquad \qquad \qquad \qquad \qquad \qquad \qquad f(a)=f'(a)=0,\quad f(b)=f'(b)=0\}.
\end{align*}
We denote it briefly by $\mathcal{S}_\eps$.

Note that for some $\Psi,\Phi,\Upsilon_1,\Upsilon_2\in C_0^\infty(-1,1)$ the function $\Psi_\eps$ converges in the sense of distributions as $\eto$ to the linear combination of the derivatives of the Dirac delta-function, which serves as a motivation for the choice of the singular perturbation $\Psi_\eps$.
If therefore the operator $\mathcal{S}_\eps$ converges (in some sense) as $\eto$ to the limit operator, then it is natural to regard this limit as the interpretation of the fourth order differential operator $\frac{d^4}{dx^4}+U+\alpha \delta'''+\beta \delta''+\gamma_1\delta'+\gamma_2\delta$.

Our purpose is to investigate the asymptotic behavior of eigenvalues $\lme$ and eigenfunctions $\ye$ as $\eto$. The perturbation $\Psi_\eps$ consists of four terms each of which has different influence on $\lme$ and $\ye$.
It is of interest to know when each term starts to have effect in asymptotic expansions.
Intuitively, we expect that the term approximating the third derivative of the Dirac delta-function has to be dominating.
We also wish to assign an operator to each collection $(\alpha,\beta,\gamma_1,\gamma_2;\Psi,\Phi,\Upsilon_1,\Upsilon_2)$.
We base the choice of the limit operator on the proximity of its eigenvalues and eigenfunctions to those of the operators $\mathcal{S}_\eps$ for sufficiently small $\eps$.

The rest of the paper is organized as follows.
Sec. 2 includes the description of the spectrum of the perturbed operators
$\mathcal{S}_\eps$.
We show that all eigenvalues are continuous functions of $\eps$ and are bounded from above.
Generally speaking, the spectrum of this family is
not bounded from below as $\eto$: for some $\Psi $ and $\alpha $ there exists a finite
number of eigenvalues converging to $-\infty$ as $\eto$.

Then Sec. 3 presents the formal asymptotic expansions for the eigenvalues and eigen\-functions of $\mathcal{S}_\eps$.
The leading terms of asymptotic expansions and the limit operators are constructed
in the section.
We introduce a spectral characteristic of the shape $\Psi $, namely,
the \textit{resonant set} $\Sigma_{\Psi }$, which is the spectrum of the eigenvalue problem \eqref{Resonant Set}.
In the case when $\alpha $
does not belong to the resonant set, the limit operator is just the direct sum of the Dirichlet operators on $(a,0)$ and $(0,b)$ respectively. In the resonant case, when $\alpha \in\Sigma_{\Psi }$, the limit operator $\mathcal{S}_{\alpha,\beta}(\Psi,\Phi)$ acts via $\mathcal{S}_{\alpha,\beta}(\Psi,\Phi)f=Lf$ on a set of functions obeying appropriate coupling conditions at the origin.

The remainders of asymptotics for eigenvalues and eigenfunctions of
$\mathcal{S}_\eps$ are constructed in Sec. 4, because we are in need of more precise asymptotics in order to prove the approximation theorems.
In this section we also analyze the effect of each singular term.
The justification and estimation of the range of validity for the approximations are
presented in Sec. 5.

\section{Spectrum of $\mathcal{S}_\eps$ and auxiliary results}
An element $f$ of $C_0^\infty(-1,1)$ is called the \textit{$\delta^{(n)}$-like shape} if
\begin{equation*}
    \eps^{-(n+1)}f(\eps^{-1}x)\to \delta^{(n)}(x)\quad\text{as}\quad\eps\to0
\end{equation*}
in $\mathcal{D}'(\mathbb{R})$-topology. Set $\langle f\rangle_k=(k!)^{-1}\int_{-\infty}^\infty \xi^kf(\xi)\,d\xi$.
It is easy to prove that a function $f\in C_0^\infty(-1,1)$ is the $\delta^{(n)}$-like shape if and only if
$\langle f\rangle_j=0$ for $j=0,\ldots,n-1$ and $\langle f\rangle_n=(-1)^n$ (
see \cite{GolovManko1} for details).
In what follows, we denote by $\mathcal{M}_n$ the set of all $\delta^{(n)}$-like shapes, i.e.
\begin{equation*}
    \mathcal{M}_n=\big\{f\in C_0^\infty(-1,1)\colon\;\;\langle f\rangle_j=0,\quad j=0,\ldots,n-1, \quad \langle f\rangle_n=(-1)^n
\big\}.
\end{equation*}

For all $\eps>0$ the spectrum of $\mathcal{S}_\eps$ is real and discrete.
Let $\{\lambda_k^\eps\}_{k=1}^\infty$ be the eigenvalues of  $\mathcal{S}_{\eps}$ enumerated in increasing order taking multiplicity into account. Suppose that $\{y_k^\eps\}_{ k=1}^\infty$ is the $L_2(a,b)$-orthonormal system of eigenfunctions.
\begin{thm}\label{Th1}
{\sl
The eigen\-values $\lambda_k^\eps$ of the operator $\mathcal{S}_\eps$ are continuous functions of $\eps\in(0,1)$. Moreover, all eigenvalues are bounded from above as $\eto$.
Let $\Psi $ change sign and $\abs{\alpha }$ be large enough; then the spectrum of $\mathcal{S}_\eps$ is unbounded from below as $\eto$, in particular, $\lambda_1^\eps\leq-c\eps^{-4}$ for some positive constant $c$.
There is at most a finite number $N^-$ of eigenvalues converging to $-\infty$ as $\eto$.}
\end{thm}
\begin{proof}
Let us consider the quadratic form
\begin{align*}
&q_{\eps}[u]=\int_a^b\bigl(|u''|^2+(U+\Psi_\eps)|u|^2\bigr)\,dx,\\
&u\in H=\{f\in W_2^2(a,b)\colon\; f(a)=f'(a)=0,\quad f(b)=f'(b)=0\},
\end{align*}
that is equicontinuous on the set of functions $u\in\mathcal{D}(q_\eps)\cap\{\|v\|=1\}$
with respect to $\eps$.
The minimax principle \cite[p. 343]{BerShu}
\begin{equation*}
\lambda_k^\eps=\inf\limits_{E_k}\,\sup\limits_{v\in
E_k,\: \norm{v}=1}q_{\eps}[v],
\end{equation*}
yields continuity of eigenvalues with respect to $\eps$.
Here $E_k$ runs over all $k$-dimensional linear subspaces of $H$, and $\|\cdot\|$ denotes the $L_2(a,b)$-norm.

Choose a subspace $E_k^*$ containing only elements vanishing in a neighborhood of the origin. Then we obtain
\begin{equation*}%\label{InEqMinMax}
\lambda_k^\eps\leq\sup\limits_{v\in E_k^*,\:
\norm{v}=1}q_{\eps}[v].
\end{equation*}
For sufficiently small $\eps$ the restriction of $q_{\eps}$ to $E_k^*$ does not depend on $\eps$. This yields bounded\-ness of the eigenvalues from above.

Suppose $\Psi $ changes sign. Let $u\in C_0^\infty(a,b)$ be a normalized function supported on an interval $[c_1,c_2]$, where $\Psi $ takes negative values. Consider the sequence $u_\eps(x)=\eps^{-1/2}u(\eps^{-1}x)$, $\norm{u_\eps}=1$, and assume that $\alpha >0$. From the minimax principle one can conclude that
\begin{multline*}
\eps^4\lambda_1^\eps\leq\eps^4 q_\eps[u_\eps]=\eps^4\int_{c_1\eps}^{c_2\eps}\Bigl(|u_\eps''|^2+U|u_\eps|^2+\Psi_\eps |u_\eps|^2\Bigr)\,dx=\\=
\int_{c_1}^{c_2}\Bigl(|u''|^2-\alpha |\Psi||u|^2\Bigr)\,d\xi+\\+
\eps\int_{c_1}^{c_2}\Big(\beta \Phi
+\gamma_1\eps\Upsilon_1+\gamma_2\eps^2\Upsilon_2\Big)|u|^2\,d\xi+
\eps^4\int_{c_1}^{c_2}U(\eps\xi)|u|^2\,d\xi.
\end{multline*}
The first integral gives a negative number for $\alpha >r$, where $$r=\int_{c_1}^{c_2}|u''(\xi)|^2\,d\xi\cdot
\biggl(\int_{c_1}^{c_2}|\Psi(\xi)||u|^2\,d\xi\biggr)^{-1},$$
while the other terms go to zero.
Thus for $\eps$ sufficiently small the estimate
$\lambda_1^\eps\leq-c\eps^{-4}$ holds with some positive $c$.
The case $\alpha <0$ may be handled in much the same way.

Let $N_\eps^-$ denote the number of negative eigenvalues of the operator $\mathcal{S}_\eps$.
Clearly, $N^-\leq \limsup_{\eto}N_\eps^-$.
It is well known \cite{BirmanSolom} that the estimate for the number of negative eigenvalues
\begin{equation*}
    N_\eps^-\leq c_0+c_1\int_a^b\abs{x}^3\abs{U(x)}\,dx+
    c_2\int_a^b\abs{x}^3\abs{\Psi_\eps(x)}\,dx,
\end{equation*}
holds, where $c_0$, $c_1$ and $c_2$ are positive constants.
The function $\Psi_\eps$ is supported on $[-\eps,\eps]$, thus
\begin{multline*}
\int_a^b |x|^3\,|\Psi_\eps(x)|\,dx=
\int_{-\eps}^{\eps}|x|^3\Big(\eps^{-4}|\alpha \Psi (\eps^{-1}x)|+
\eps^{-3}|\beta \Phi (\eps^{-1}x)|+\\+\eps^{-2}|\gamma_1\Upsilon_1(\eps^{-1}x)|+
\eps^{-1}|\gamma_2\Upsilon_2(\eps^{-1}x)|\Big)\,dx=
\abs{\alpha }\int_{-1}^{1} |\xi|^3\,|\Psi (\xi)|\,d\xi+\\+
\eps\int_{-1}^{1}|\xi|^3\Big(|\beta \Phi (\xi)|+
\eps|\gamma_1\Upsilon_1(\xi)|+\eps^{2}|\gamma_2\Upsilon_2(\xi)|
\Big)\,d\xi\leq|\alpha |c(\Psi )+1
\end{multline*}
for small $\eps>0$.
From what has already been proved it follows that
$$N^-\leq c_0+c_1(U)+c_2(\Psi )|\alpha |$$
for some positive constants $c_0$, $c_1(U)$, $c_2(\Psi )$.
\end{proof}

Therefore the spectrum of $\mathcal{S}_\eps$ consists of two parts: the set of eigenvalues tending to $-\infty$ as $\eto$, and the set of all bounded eigenvalues as $\eto$.

\section{Asymptotics of eigenvalues and eigenfunctions of $\mathcal{S}_\eps$ and the limit operator}\label{SecMT}
Fix an eigenvalue $\lambda_k^\eps$ of the problem \eqref{MainProblem}
with $k>N^-$. We write it $\lambda^\eps$ for short. Let $\ye$ be the corresponding eigenfunction.
The asymptotic expansions of $\lme$ are represented by
\begin{align}\label{ExpanEValue}
    \lambda^\eps\sim&\, \lambda+\eps\lambda_1+\eps^2\lambda_2+\ldots,\\
\intertext{and we postulate two-scale expansions for the eigenfunction}
\label{ExpanEFunctBigSize}
    \ye(x)\sim&\, v(x)+\eps v_1(x)+\eps^2v_2(x)+\ldots\quad \text{for}\quad x\in(a,-\eps)\cup(\eps,b),\\
\label{ExpanEFunctSmallSize}
    \ye(x)\sim&\, \eps w(\eps^{-1}x)+\eps^2w_1(\eps^{-1}x)+\ldots\quad \text{for}\quad x\in(-\eps,\eps).
\end{align}
Here all functions $v$, $v_k$ are defined for $x\in(a,0)\cup(0,b)$, and $w$, $w_k$ are defined for $\xi\in(-1,1)$. Assume that $v$ is different from zero. Series \eqref{ExpanEFunctBigSize}, \eqref{ExpanEFunctSmallSize} satisfy the coupling conditions
\begin{equation}\label{FittConditions}
    \Big[\ye^{(j)}\Big]_{x=-\eps}=0,\qquad\Big[\ye^{(j)}\Big]_{x=\eps}=0,\qquad j=0,\ldots ,3,
\end{equation}
where by $[f]_{x=a}$ we denote the jump of $f$ at a point $a$.

We substitute series \eqref{ExpanEValue}, \eqref{ExpanEFunctBigSize} into the equation and the boundary conditions \eqref{MainProblem} and derive
\begin{align}\label{vEquation}
  &Lv=\lambda v, &&v(a)=v'(a)=0,\quad v(b)=v'(b)=0,\\
  \label{v1Equation}
  &Lv_1=\lambda v_1+\lambda_1v, && v_1(a)=v_1'(a)=0,\quad v_1(b)=v_1'(b)=0,\\  \label{v2Equation}
  &Lv_2=\lambda v_2+\lambda_1v_1+\lambda_2v, && v_2(a)=v_2'(a)=0,\quad v_2(b)=v_2'(b)=0,
\end{align}
where all equations hold on $(a,0)\cup(0,b)$.
We set $\xi=\eps^{-1}x$. After substituting \eqref{ExpanEValue}, \eqref{ExpanEFunctSmallSize} into the equation \eqref{MainProblem},
one obtains the following equations on $(-1,1)$
\begin{align}
&w^{(4)}+\alpha \Psi  w=0,\\
&w_1^{(4)}+\alpha \Psi  w_1=-\beta \Phi  w,\label{w1Equation}\\
&w_2^{(4)}+\alpha \Psi  w_2=-\beta \Phi  w_1-\gamma_1\Upsilon_1 w,\label{w2Equation}\\
&w_3^{(4)}+\alpha \Psi  w_3=-\beta \Phi  w_2-\gamma_1\Upsilon_1 w_1-\gamma_2\Upsilon_2 w.\label{w3Equation}
\end{align}
Substituting \eqref{ExpanEFunctBigSize}, \eqref{ExpanEFunctSmallSize} into the coupling conditions \eqref{FittConditions}, we can assert that
\begin{equation*}
v^{(j)}(\pm\eps)+\eps \,v_1^{(j)}(\pm\eps)+\eps^2 \,v_2^{(j)}(\pm\eps)+\cdots \sim
        \eps^{1-j}\,w^{(j)}(\pm1)+\eps^{2-j}\, w_1^{(j)}(\pm1)+\cdots
\end{equation*}
for $j=0,\ldots,3$.
We can now expand $v_k^{(j)}$ into the formal Taylor series about $x=\pm0$.
Then we conclude that
\begin{align}\label{v(0)}
v(-0)=0,&\qquad v(+0)=0,\\
w''(\pm1)=0,&\qquad w'''(\pm1)=0,\\
\label{AdCondV}
v'(-0)=w'(-1),&\qquad v'(+0)=w'(1),\\
\label{AdCondVPrime}
v_1(-0)- v '(-0)=w(-1),& \qquad v_1(+0)+ v '(+0)=w(1),\\
\label{AdCondVPrimeVW}
v'_1(-0)- v''(-0)=w'_1(-1),& \qquad v'_1(+0)+ v''(+0)=w'_1(1),\\
\label{AdCondV2W}
v_2(-0)-v_1'(-0)+\textstyle\frac{1}{2}v''(-0)=w_1(-1),& \qquad v_2(+0)+v_1'(+0)+\textstyle\frac{1}{2}v''(+0)=w_1(1),\\
\label{AdCondVPrimeV2W}
v'_2(-0)-v_1''(-0)+\textstyle\frac{1}{2}v'''(-0)=w'_2(-1),& \qquad v'_2(+0)+v_1''(+0)+\textstyle\frac{1}{2}v'''(+0)=w'_2(1),\\
\label{AdCondW1Primes}
w''_1(\pm1)=v''(\pm0),&\qquad w'''_1(\pm1)=0,\\
\label{AdCondW2Primes2}
w''_2(\pm1)=v''_1(\pm0)\pm v'''(\pm0),&\qquad w'''_2(\pm1)=v'''(\pm0),\\
\label{AdCondW2Primes3}
w''_3(\pm1)=v''_2(\pm0)\pm v'''_1(\pm0)+\textstyle\frac{1}{2}v^{(4)}(\pm0),&\qquad w'''_3(\pm1)=v'''_1(\pm0)\pm v^{(4)}(\pm0).
\end{align}
It follows that $v$ satisfies the equation and the boundary conditions \eqref{vEquation},
and furthermore $v(0)=0$. The function $w$ is a solution to the problem
\begin{equation}\label{w0Problem}
w^{(4)}+\alpha \Psi  w=0,\quad\xi\in(-1,1),\qquad
w''(\pm1)=0,\quad w'''(\pm1)=0.
\end{equation}
Moreover these functions are related by the coupling conditions \eqref{AdCondV}.
The problem \eqref{w0Problem} is decisive in our next consideration, because it contains information about the singular perturbation. The first and primary question is whether there exists its nontrivial solution.

\subsection{Resonant set}
The problem \eqref{w0Problem} can be regarded as a spectral problem with the spectral parameter $\alpha $.
We note that in the generic case, the function $\Psi $ is sign-changing. It is of interest to investigate spectral properties of this problem. We will also introduce the spectral characteristic of the shape $\Psi $.

We introduce the operator $\mathcal{T}_\Psi=\frac{1}{\Psi (\xi)}\frac{d^4}{d\xi^4}$ with the domain
$$
\mathcal{D}(\mathcal{T}_{\Psi})=\{f\in L^2_{|\Psi |}(-1,1)\mid \Psi ^{-1}f^{(4)}\in L^2_{|\Psi |}(-1,1),\quad f''(\pm1)=f'''(\pm1)=0
\}.
$$
The problem \eqref{w0Problem} is equivalent to the spectral equation $\mathcal{T}_{\Psi}w=-\alpha w$.
\begin{thm}\label{ThJproperties}
{\sl
Given $\Psi \in C_0^\infty(\mathbb{R})$, with $\mathrm{supp}\,\Psi=[-1,1]$, the spectrum of the operator $\mathcal{T}_\Psi$  is real and discrete.
Suppose $\Psi$ changes sign; then the spectrum of $\mathcal{T}_\Psi$ has two accumulation points $-\infty$ and $+\infty$.
}\end{thm}
\begin{proof}
Since the case where $\Psi$ keeps sign is much simpler and can be handled within the standard Hilbert space theory, we assume that $\Psi$ changes sign and apply the Krein space theory to investigate the spectrum of $\mathcal{T}_\Psi$.
Let $\mathcal{L}$ be the weighted $L_2$-space with the scalar product $(f,g)=\int_{-1}^1|\Psi |f\bar{g}\,d\xi$,
and let us define by $[f,g]=\int_{-1}^1\Psi  f\bar{g}\,d\xi$ the indefinite metric in $\mathcal{L}$. Then the pair $(\mathcal{L},[\cdot,\cdot])$ is called a Krein space \cite[ch. 1]{IA}.

In this Krein space there exists the fundamental symmetry $Jf=\mathrm{sgn}\,{\Psi }f$ such that $[f,g]=(Jf,g)$ for all $f,g\in\mathcal{L}$.
An operator $T$ is $J$-\textit{selfadjoint} if $JT$ is selfadjoint in $L^2_{|\Psi |}(-1,1)$. An operator $T$ is said to be $J$-\textit{nonnegative} if $[Tf,f]\geq0$ for all $f\in\mathcal{D}(T)$.

For each $\Psi \in C_0^\infty(-1,1)$ the operator $J\mathcal{T}_{\Psi}$ is selfadjoint, and so $\mathcal{T}_{\Psi}$ is $J$-selfadjoint.
Next, for all $f\in \mathcal{D}(\mathcal{T}_{\Psi})$ one obtains
$$
[\mathcal{T}_{\Psi} f,f]=\int_{-1}^1f^{(4)} \overline{f}\,d\xi=
\int_{-1}^1\abs{f''}^2\,d\xi\geq0.
$$
Thus $\mathcal{T}_{\Psi}$ is $J$-nonnegative.
Any $J$-selfadjoint and $J$-nonnegative operator with a nonempty resolvent set has real spectrum \cite[p. 138]{IA}. Let us show that the resolvent set of $\mathcal{T}_{\Psi}$ is nonempty. The homogenous problem
\begin{equation}
\label{HomogProb}
g^{(4)}+i\Psi  g=0,\quad\xi\in(-1,1),\qquad
g''(\pm1)=0,\quad g'''(\pm1)=0
\end{equation}
has a trivial solution only. Indeed, each solution satisfies the equality
$$
\int_{-1}^1\abs{g''}^2d\xi+i\int_{-1}^1\Psi \abs{g}^2d\xi=0.
$$
Since $\Psi $ is real-valued, it follows that $g$ is a linear function. Obviously, only zero function can be a solution of \eqref{HomogProb}.
Hence the nonhomogeneous problem $g^{(4)}+i\Psi  g=h$, $g''(\pm1)=0$, $ g'''(\pm1)=0$ admits a unique solution for arbitrary $h\in L_2(-1,1)$ \cite[p. 39]{Na}. Note that $\Psi f$ belongs to $L_2(-1,1)$  for each $f\in \mathcal{L}$, since $\norm{\Psi f}_{L_2(-1,1)}\leq\max_\mathbb{R}|\Psi |^{1/2}\cdot\norm{f}_\mathcal{L}$. Then the equation
$\mathcal{T}_{\Psi}g+i g=f$ is equivalent to
the nonhomogeneous problem $g^{(4)}+i\Psi  g=\Psi f$, $g''(\pm1)=0$, $ g'''(\pm1)=0$ and  admits a unique solution for each $f\in \mathcal{L}$.
Therefore  $-i$ belongs to the resolvent set. Since the resolvent set of  $\mathcal{T}_{\Psi}$ is nonempty, the spectrum of $\mathcal{T}_{\Psi}$ is real.

We shall prove that the resolvent $R_\mu(\mathcal{T}_{\Psi})$ of the operator $\mathcal{T}_{\Psi}$ is compact. The operator $R_\mu(\mathcal{T}_{\Psi})$ acts from the space $\mathcal{L}$ into $\mathcal{D}(\mathcal{T}_{\Psi})$, and for each $f\in\mathcal{L}$ solves the equation
\begin{equation*}
    g^{(4)}-\mu\Psi g=\Psi f,\qquad g\in \mathcal{D}(\mathcal{T}_{\Psi}).
\end{equation*}
As far as the right-hand side $\Psi f$ belongs to $L_2(-1,1)$, it follows that the solution $g$ is an element of $W_2^4(-1,1)$.
The space $\mathcal{D}(\mathcal{T}_{\Psi})$ is a Banach space with the graph norm.
The sequence of continuous embeddings
$\mathcal{D}(\mathcal{T}_{\Psi})\subset W_2^4(-1,1)\subset L_2(-1,1)\subset
\mathcal{L}$ yields the compactness of the resolvent, since
$W_2^4(-1,1)\subset L_2(-1,1)$ is the compact embedding.
As a consequence we have $\sigma(\mathcal{T}_{\Psi})=\sigma_p(\mathcal{T}_{\Psi})$.

Since $\Psi $ changes sign, the spectrum $\sigma(\mathcal{T}_{\Psi})$ is unbounded in both directions \cite{CurJDE}.
\end{proof}

We introduce the set $\Sigma_{\Psi }=\{\alpha \in\mathbb{R}\colon-\alpha \in\sigma(\mathcal{T}_{\Psi})\}$, which is the spectrum of the problem \eqref{w0Problem}. We call $\Sigma_{\Psi }$ the \textit{resonant set} of the shape $\Psi $.
When $\alpha\in\Sigma_\Psi$, suppose that \eqref{MainAssumption} holds
(the case of {\it nondegenerate resonance}).
In this paper we assume that only the nondegenerate resonance is possible, namely,
if $\alpha $ belongs to the resonant set, then both conditions \eqref{MainAssumption} hold.

\subsection{The limit operator}
Let us continue to construct the asymptotics. We distinguish two different cases and
start with the assumption $\alpha \notin\Sigma_{\Psi }$. Then the problem
\eqref{w0Problem} admits a trivial solution $w=0$ only. That $v'(0)=0$ follows from the coupling conditions \eqref{AdCondV}. We conclude from \eqref{vEquation} that $v$ is a solution to the problem
\begin{equation}\label{LimitProblemw=0}
\begin{cases}
 Lv=\lambda v,\qquad x\in(a,0)\cup(0,b),\\
 v(a)=v'(a)=0,\quad v(0)=v'(0)=0,\quad v(b)=v'(b)=0.
 \end{cases}
\end{equation}
Let us introduce the operators
\begin{align*}
    &S_-f=Lf,\quad \mathcal{D}(S_-)=\{f\in W_2^4(a,0)\colon\;\; f(a)=f'(a)=0,\quad f(0)=f'(0)=0\},\\
    &S_+f=Lf,\quad \mathcal{D}(S_+)=\{f\in W_2^4(0,b)\colon\;\; f(0)=f'(0)=0,\quad f(b)=f'(b)=0\}.
\end{align*}
The operator $S_-\oplus S_+$ is associated with the problem \eqref{LimitProblemw=0}.
Therefore in the nonresonant case, when $\alpha \notin\Sigma_{\Psi }$, we can define the limit operator as $S_-\oplus S_+$.

Let us now suppose that $\alpha $ belongs to the resonant set $\Sigma_{\Psi }$. Recalling \eqref{MainAssumption}, we deduce that the quotient
\begin{equation*}
    \theta_{\Psi }(\alpha )=\frac{w'_{\alpha }(1)}{w'_{\alpha }(-1)}
\end{equation*}
is well defined and does not depend on the choice of an eigenfunction.
Clearly, $w=cw_{\alpha }(\xi)$, where $c$ is a constant.
We conclude from \eqref{AdCondV} that $v'(-0)=cw'_{\alpha }(-1)$, $v'(+0)=cw'_\alpha(1)$, hence that
\begin{equation}\label{ThetaVCond1}
    v'(+0)-\theta_{\Psi }(\alpha ) v'(-0)=0,
\end{equation}
and also that $c=\frac{v'(-0)}{w'_{\alpha }(-1)}$. According to \eqref{w1Equation}, \eqref{AdCondW1Primes} the next term $w_1$ of series \eqref{ExpanEFunctSmallSize} can be found by solving the problem
\begin{equation}\label{w1Problem}
\begin{aligned}
&w_1^{(4)}+\alpha \Psi  w_1=-\beta \textstyle\frac{v'(-0)}{w'_{\alpha }(-1)}\Phi  w_{\alpha },\qquad\xi\in(-1,1),\\
&w''_1(-1)=v''(-0),\quad w_1'''(-1)=0,\quad
w''_1(1)=v''(+0),\quad w_1'''(1)=0.
\end{aligned}
\end{equation}
Because $\alpha $ is an eigenvalue of \eqref{w0Problem}, the problem admits a solution if and only if
\begin{equation}\label{ThetaVCond2}
    \theta_{\Psi }(\alpha )v''(+0)-v''(-0)= \beta  v'(-0)
    \int_{-1}^1\Phi (\xi)\textstyle\left(\frac{w_{\alpha }(\xi)}{w'_{\alpha }(-1)} \right)^2 \,d\xi.
\end{equation}
To derive this solvability condition, we multiply the equation by $w_{\alpha }$ and integrate by parts.
Let us define a functional on $C_0^\infty(-1,1)$ by
$\vartheta_\Phi[f]=\int_{-1}^1\Phi f^2\,d\xi$.
Collecting \eqref{vEquation}, \eqref{ThetaVCond1} and \eqref{ThetaVCond2}
we deduce that $v$ must be an eigenfunction of the problem
\begin{equation}\label{LimitProblemV}
\begin{cases}
    Lv=\lambda v,\qquad x\in(a,0)\cup(0,b),\\
    v(a)=v'(a)=0,\quad  v(b)=v'(b)=0,\\
    v(0)=0,\quad v'(+0)-\theta_{\Psi }(\alpha ) v'(-0)=0,\\
    \theta_{\Psi }(\alpha ) v''(+0)-v''(-0)-\beta \vartheta_\Phi\big[w_{\alpha }/w'_{\alpha }(-1)\big] v'(-0)=0.
\end{cases}
  \end{equation}
Consequently, the operator $S(\alpha,\beta;\Psi,\Phi)=\frac{d^4}{dx^4}+U(x)$ with the domain
\begin{multline}
\label{LimitOperator}
    \mathcal{D}(S(\alpha,\beta;\Psi,\Phi))=
    \Big\{f\in W_2^4\big((a,0)\cup(0,b)\big)\colon\;
    f(a)=f'(a)=0,\\ f(b)=f'(b)=0,\quad
    f(0)=0,\quad
    f'(+0)-\theta_{\Psi }(\alpha ) f'(-0)=0,\\
    \theta_{\Psi }(\alpha ) f''(+0)-f''(-0)-\beta \vartheta_\Phi\big[w_{\alpha }/w'_{\alpha }(-1)\big] f'(-0)=0
    \Big\}
\end{multline}
is associated with the problem \eqref{LimitProblemV}. Combining resonant case and nonresonant one, gives us the limit operator
\begin{equation*}
\mathcal{S}_{\alpha,\beta}(\Psi,\Phi)=
    \begin{cases}
    S_-\oplus S_+,& \alpha\notin\Sigma_\Psi,\\
    S(\alpha,\beta;\Psi,\Phi),& \alpha\in\Sigma_\Psi.
    \end{cases}
\end{equation*}
Recall that we consider only those $\alpha$ from the resonant set, which satisfy
assumptions \eqref{MainAssumption}.

\section{Asymptotic expansions of eigenvalues and eigenfunctions of $\mathcal{S}_\eps$ : correctors}\label{SubsecCorr}
In order to justify the closeness of eigenvalues and eigenfunctions of operators $\mathcal{S}_\eps$ and $\mathcal{S}_{\alpha,\beta}(\Psi,\Phi)$ we must derive next terms of series \eqref{ExpanEValue}--\eqref{ExpanEFunctSmallSize}.
Clearly, the  construction of correctors depends on the multiplicity of $\lambda$.
Let $\lambda$ be a simple eigenvalue of $\mathcal{S}_{\alpha,\beta}(\Psi,\Phi)$ with the eigenfunction $v$ being normalized in $L_2(a,b)$.
\subsection{Asymptotics in the nonresonant case}
In this subsection we assume that $\alpha $ does not belong to the resonant set $\Sigma_{\Psi }$. Then $w=0$ and $\sigma(\mathcal{S}_{\alpha,\beta}(\Psi,\Phi))=\sigma(S_-)\cup\sigma(S_+)$.
If $\lambda$ is a simple eigenvalue of $\mathcal{S}_{\alpha,\beta}(\Psi,\Phi)$, then $\lambda$ is a simple eigenvalue of $S_-$ or $S_+$.
Without loss of generality we may assume $\lambda\in\sigma(S_+)$, and thus $v$ vanishes on $(a,0)$. Employing \eqref{w1Equation}, \eqref{AdCondW1Primes} gives us the problem
\begin{align*}
&w_1^{(4)}+\alpha \Psi  w_1=0,\qquad\xi\in(-1,1),\\
&w_1''(-1)=w_1'''(-1)=0,\quad w_1''(1)=v''(+0),\quad w'''_1(1)=0,
\end{align*}
which admits a unique solution, since $\alpha $ does not belong to the spectrum of \eqref{w0Problem}.
In light of \eqref{v1Equation}, \eqref{AdCondVPrime} the function $v_1$ can be found by solving problems
\begin{equation}\label{ProblemV1-}
\begin{aligned}
  &\begin{cases}
  Lv_1=\lambda v_1, \quad x\in(a,0),\\
  v_1(a)=v'_1(a)=0,\\
  v_1(-0)=0,\quad v_1'(-0)=w_1'(-1),
  \end{cases}
\\
  &\begin{cases}
  Lv_1=\lambda v_1+\lambda_1 v, \quad x\in(0,b),\\
  v_1(+0)=-v '(+0),\quad v_1'(+0)=w_1'(1)-v''(+0),\\
v_1(b)=v'_1(b)=0.
  \end{cases}
\end{aligned}
\end{equation}
on $(a,0)$ and $(0,b)$ respectively.
Of course, the first of these problems has a unique solution, since $\lambda\notin\sigma(S_-)$.
Note that in the generic case  the second problem has no solution. But  we
can ensure the existence of a solution by choosing the free parameter $\lambda_1$. Indeed, applying the Fredholm alternative we conclude that the second problem \eqref{ProblemV1-} admits a solution if and only if
$$\lambda_1=v ''(+0)(v''(+0)-w'_1(1))-v'(+0)v'''(+0).$$
To derive this we multiply the equation by the eigenfunction and integrate by parts.
The last equality is simultaneously a formula for the corrector $\lambda_1$ in the asymptotic expansions of the eigenvalue.
Clearly, the solution $v_1$ is defined up to the term $cv$. To fix it we subordinate the solution to the condition $\int_0^bv v_1\,dx=0$.

Combining \eqref{w2Equation} with \eqref{AdCondW2Primes2} and recalling $\alpha \notin\Sigma_{\Psi }$, we deduce the problem
\begin{align*}%\label{ProblemW2}
    &w_2^{(4)}+\alpha \Psi  w_2=-\beta \Phi  w_1 ,\qquad \xi\in(-1,1),\\
&w''_2(-1)=v''_1(-0),\quad w'''_2(-1)=0,\\
&w''_2(1)=v''_1(+0)+ v'''(+0),\quad w'''_2(1)=v'''(+0),
\end{align*}
which gives us the corrector $w_2$.
We employ \eqref{v2Equation}, \eqref{AdCondV2W}, \eqref{AdCondVPrimeV2W} to find
\begin{equation}\label{ProblemV2-}
\begin{aligned}
  &\begin{cases}
  Lv_2=\lambda v_2+\lambda_1 v_1, \quad x\in(a,0),\\
  v_2(a)=0,\quad v'_2(a)=0,\\
  v_2(-0)=w_1(1)+v_1'(-0),\\
  v_2'(-0)=w_2'(-1)+v_1''(-0),
  \end{cases}
\\
  &\begin{cases}
  Lv_2=\lambda v_2+\lambda_1 v_1+\lambda_2 v, \quad x\in(0,b),\\
  v_2(+0)=w_1(1)-v_1 '(+0)-\frac{1}{2}v''(+0),\\
  v_2'(+0)=w_2'(1)-v_1''(+0)-\frac{1}{2}v'''(+0),\\
v_2(b)=v'_2(b)=0.
  \end{cases}
\end{aligned}
\end{equation}
As before we deduce that the first of these problems has a unique solution, and the second one admits a solution if and only if
\begin{equation*}
    \lambda_2=v'''(+0)\bigl(w_1(1)-v_1'(+0)-\textstyle\frac{1}{2}v''(+0)\bigr)-
    v''(+0)\bigl(w_2'(1)-v_1''(+0)-\frac{1}{2}v'''(+0)\bigr).
    \end{equation*}
For the sake of definiteness, the solution is subject to the additional condition $\int_0^bvv_2\,dx=0$.
By using \eqref{w3Equation}, \eqref{AdCondW2Primes3} one obtains
\begin{align*}%\label{ProblemW3}
&w_3^{(4)}+\alpha \Psi  w_3=-\beta \Phi  w_2-\gamma_1\Upsilon_1 w_1 ,\qquad \xi\in(-1,1),\\
&w''_3(-1)=v_2''(-0)- v_1'''(-0),\quad w'''_3(-1)=v'''_1(-0),\\
&w''_3(1)=v_2''(+0)+ v_1'''(+0)+\textstyle\frac{1}{2}v^{(4)}(+0),\quad w'''_3(1)=v'''_1(+0)+ v^{(4)}(+0).
\end{align*}
Reasoning as before, from this problem we get $w_3$.

Let us introduce the notations
\begin{equation}\label{Approx}
\begin{gathered}
    \Lambda_\eps=\lambda +\eps\lambda_1+\eps^2\lambda_2, \\
    Y_\eps(x)=
    \begin{cases}
    v (x)+\eps v_1(x)+\eps^2 v_2(x),  &x\in(a,-\eps)\cup(\eps,b),\\
    \eps^2 w_1(\eps^{-1}x)+\eps^3 w_2(\eps^{-1}x)+\eps^4 w_3(\eps^{-1}x),&x\in(-\eps,\eps)
    \end{cases}
\end{gathered}
\end{equation}
for the constructed approximations of eigenvalues and eigenfunctions.

\subsection{Asymptotics under resonance}

Now we assume that $\alpha $ belongs to the resonant set $\Sigma_{\Psi }$ and that
$\lambda$ is an eigenvalue of the operator $S(\alpha,\beta,\Psi,\Phi)$.
Let $w_{\alpha }$ be an eigenfunction of \eqref{w0Problem} such that
$w'_{\alpha }(-1)=1$. Clearly, $\theta_{\Psi }(\alpha )=w'_{\alpha }(1)$.

Since \eqref{ThetaVCond2} holds, the problem \eqref{w1Problem} admits a solution. This solution can be represented as $w_1=w_1^*+c_1w_{\alpha }$, the function $w_1^*$ being a partial solution of the problem fixed by the condition $\frac{dw_1^*}{d\xi}(-1)=0$, and the constant $c_1$ is to be chosen later.

We next construct the corrector $v_1$.
The function $v_1$ satisfies the equation \eqref{v1Equation} outside the origin  and \eqref{AdCondW1Primes} yields
\begin{equation}\label{V1D1MatchCond}
    v'_1(+0)-\theta_{\Psi }(\alpha ) v'_1(-0)=G_1,
\end{equation}
where $G_1=w'_1(1)-\theta_{\Psi }(\alpha )w'_1(-1)-v''(+0)-\theta_{\Psi }(\alpha )v''(-0)$.
Although $w_1$ is not uniquely chosen, the constant $G_1$ is well defined.
In fact,
\begin{multline*}
w'_1(1)-\theta_{\Psi }(\alpha ) w'_1(-1)=\biggl(\frac{dw_1^*}{d\xi}(1)-
\theta_{\Psi }(\alpha )\frac{dw_1^*}{d\xi}(-1)\biggr)
    +\\+c_1\bigl(w'_{\alpha }(1)-\theta_{\Psi }(\alpha )w'_{\alpha }(-1)\bigr)=
    \frac{dw_1^*}{d\xi}(1).
\end{multline*}
From \eqref{w2Equation} and \eqref{AdCondW2Primes2} it follows that the corrector $w_2$ must solve the problem
\begin{align}\label{ProblemW2RC}
\begin{aligned}
&w_2^{(4)}+\alpha \Psi  w_2=-\beta \Phi w_1-\gamma_1\Upsilon_1w,\qquad \xi\in(-1,1),\\
&w''_2(-1)=v''_1(-0)- v'''(-0),\quad w'''_2(-1)=v'''(-0),\\
&w''_2(1)=v''_1(+0)+ v'''(+0),\quad w'''_2(1)=v'''(+0).
\end{aligned}
\end{align}
From the first condition in \eqref{AdCondVPrimeVW} we deduce  $c_1=v_1'(-0)-v''(-0)$.
Set $\vartheta_{\Upsilon_i}[f]=\int_{-1}^1\Upsilon_if\,d\xi$ for $f\in C_0^\infty(-1,1)$ and $i=1,2$.
Thus the solvability condition of the above problem can be written as
\begin{equation}\label{SolvabilityW2}
    \theta_{\Psi }(\alpha ) v_1''(+0)-v_1''(-0)-\beta \vartheta_\Psi[w_{\alpha }]v_1'(-0)=H_1
\end{equation}
which is due to the Fredholm alternative. Here
\begin{multline*}
H_1=\big(w_{\alpha }(1)-\theta_{\Psi }(\alpha )\big)v'''(+0)-
\big(w_{\alpha }(-1)+1\big)v'''(-0)+  \\+
\beta \big(\vartheta_\Psi[\sqrt{w_1^*w_\alpha} ]-\vartheta_\Psi[w_{\alpha }]v''(-0)\big)+
\gamma_1\vartheta_{\Upsilon_1}[w^2_{\alpha }]v'(-0).
\end{multline*}
From \eqref{AdCondVPrime} we have $v_1(\pm0)=v'(-0)w_{\alpha }(\pm1)\mp v'(\pm0)$.
Combining these identities along with \eqref{v1Equation}, \eqref{V1D1MatchCond} and \eqref{SolvabilityW2} we conclude that $v_1$ solves the problem
\begin{equation}\label{ProblemV1}
\begin{cases}
    Lv_1=\lambda v_1+\lambda_1v,\qquad x\in(a,0)\cup(0,b),\\
    v_1(a)=v'_1(a)=0,\quad v_1(b)=v'_1(b)=0,\\
    v_1(-0)=v'(-0)w_{\alpha }(-1)+v'(-0),\quad v_1(+0)=v'(-0)w_{\alpha }(1)- v'(+0),\\
    v'_1(+0)-\theta_{\Psi }(\alpha ) v'_1(-0)=G_1,\quad
    \theta_{\Psi }(\alpha ) v_1''(+0)-v_1''(-0)-\beta \vartheta_\Psi[w_{\alpha }]v_1'(-0)=H_1.
\end{cases}
\end{equation}
The free parameter $\lambda_1$ in the right-hand side of equation \eqref{ProblemV2} enables us to solve the problem.
In view of Fredholm's alternative, \eqref{ProblemV2} admits a solution if and only if
\begin{multline*}
\lambda_1=H_1v'(-0)-G_1v''(+0)-\big(v'(-0)w_{\alpha }(-1)+v'(-0)\big)v'''(-0)+\\+
\big(v'(-0)w_{\alpha }(1)- v'(+0)\big)v'''(+0).\end{multline*}
For the sake of definiteness, the solution is subject to the additional condition
$\int_a^bv v_1\,dx=0$.

Given $v_1$, we may compute the constant $c_1$. A trivial verification shows that the second condition in \eqref{AdCondVPrimeVW} holds.

The condition \eqref{SolvabilityW2} enables one to solve the problem \eqref{ProblemW2RC}. A solution of this problem has the form $w_2=w_2^*+c_2w_{\alpha }$, where $w_2^*$ solves \eqref{ProblemW2RC} and satisfies $\frac{dw_2^*}{d\xi}(-1)=0$. The constant $c_2$ will be chosen later.

Following as before we shall similarly find correctors $w_3$, $v_2$ and $\lambda_2$.
The function $v_2$ satisfies the equation and boundary conditions \eqref{v2Equation}, and
\begin{equation}\label{G2}
    v'_2(+0)-\theta_{\Psi }(\alpha )v'_2(-0)=G_2
\end{equation}
by \eqref{AdCondVPrimeV2W}. Here $G_2=\frac{dw_2^*}{d\xi}(1)-\theta_{\Psi }(\alpha )\big(v_1''(-0)-\frac{1}{2}v'''(-0)\big)-
v_1''(+0)-\frac{1}{2}v'''(+0)$.
Next we employ \eqref{w3Equation} and \eqref{AdCondW2Primes3} to obtain the problem
\begin{align}\label{LastW3Problem}
\begin{aligned}
&w_3^{(4)}+\alpha \Psi w_3=-\beta \Phi w_2-\gamma_1\Upsilon_1w_1-\gamma_2\Upsilon_2w,\qquad \xi\in(-1,1),\\
&w_3''(-1)=v_2''(-0)-v_1'''(-0)+\textstyle\frac{1}{2}v^{(4)}(-0),\quad
w_3'''(-1)=v_1'''(-0)-v^{(4)}(-0),\\
&w_3''(1)=v_2''(+0)+v_1'''(+0)+\textstyle\frac{1}{2}v^{(4)}(+0),\quad
w_3'''(1)=v_1'''(+0)+v^{(4)}(+0).
\end{aligned}
\end{align}
On applying \eqref{AdCondVPrimeV2W}, one obtains $c_2=v_2'(-0)-v_1''(-0)+\textstyle\frac{1}{2}v'''(-0)$. Therefore we may write the solvability condition of this problem in the form
\begin{equation}\label{V2D2MatchCond}
    \theta_{\Psi }(\alpha )v_2''(+0)-v_2''(-0)-\beta \vartheta_\Psi[w_{\alpha }]v_2'(-0)=
    H_2
\end{equation}
with
\begin{multline*}
    H_2=v_1'''(+0)\big\{w_{\alpha }(1)-\theta_{\Psi }(\alpha )\big\}-
    v_1'''(-0)\big\{w_{\alpha }(-1)+1\big\}+\\+
    v^{(4)}(+0)\big\{w_{\alpha }(1)-\textstyle\frac{1}{2}\theta_{\Psi }(\alpha )\big\}+
    v^{(4)}(-0)\big\{w_{\alpha }(-1)+\textstyle\frac{1}{2}\big\}+\\+
    \beta \bigl(\vartheta_\Psi[\sqrt{w_2^*w_{\alpha }}]-\bigl\{v_1''(-0)-\textstyle\frac{1}{2}v'''(-0) \bigr\}\vartheta_\Psi[w_{\alpha }]\bigr)+\gamma_1\vartheta_{\Upsilon_1}[w_1w_{\alpha }]+
    \gamma_2\vartheta_{\Upsilon_2}[ww_{\alpha }].
\end{multline*}
From \eqref{AdCondV2W}
we find $v_2(\pm0)=F_2^\pm$, where
$$F_2^\pm=w_1^*(\pm1)+\big(v_1'(-0)-v''(-0)\big)w_{\alpha }(\pm1)\mp v_1'(\pm0)-\textstyle\frac{1}{2}v''(\pm0).$$
In view of \eqref{v2Equation}, \eqref{G2} and
\eqref{V2D2MatchCond} it follows that $v_2$ is a solution to the problem
\begin{equation}\label{ProblemV2}
\begin{cases}
    Lv_2=\lambda v_2+\lambda_1v_1+\lambda_2v,\qquad x\in(a,0)\cup(0,b),\\
    v_2(a)=v'_2(a)=0,\quad v_2(b)=v'_2(b)=0,\\
    v_2(-0)=F_2^-,\quad v_2(+0)=F_2^+,\quad v'_2(+0)-\theta_{\Psi }(\alpha ) v'_2(-0)=G_2,\\
        \theta_{\Psi }(\alpha ) v_2''(+0)-v_2''(-0)-\beta \vartheta_\Psi[w_{\alpha }]v_2'(-0)=H_2.
\end{cases}
\end{equation}
The problem admits a solution if and only if $$\lambda_2=H_2v'(-0)-G_2v''(+0)-F_2^-v'''(-0)+F_2^+v'''(+0).$$ This solution is defined up to the term $cv$. To eliminate this ambiguity we additionally demand that the condition $\int_a^bv v_2\,dx=0$ holds.

Summing up, one obtains the following approximations for the eigenvalue and eigenfunction of the perturbed problem in the resonant case:
\begin{equation}\label{Approx1}
\begin{gathered}
    \Lambda_\varepsilon=\;\lambda +\varepsilon\lambda_1+\varepsilon^2\lambda_2, \\ Y_\varepsilon(x)=
\begin{cases}
v (x)+\varepsilon v_1(x)+\varepsilon^2 v_2(x),  &x\in(a,-\eps)\cup(\eps,b),\\
\eps w(\varepsilon^{-1}x) + \varepsilon^2
w_1(\varepsilon^{-1}x)+\varepsilon^3 w_2(\varepsilon^{-1}x)+\eps^4w_3(\varepsilon^{-1}x),&x\in(-\eps,\eps).
\end{cases}
\end{gathered}
\end{equation}
Here $w_3$ is an arbitrary solution of \eqref{LastW3Problem}.
The choice of  $c_3$ in the representation $w_3=w_3^*+c_3w_\alpha$ is not important since we do not look for the corrector $v_3$.

\section{Justification of asymptotic expansions}
As shown in Theorem~\ref{Th1}, for every regularization $\Psi_\eps(x)$
there is at most a finite number
$N^->0$ of eigenvalues $\lme_k$, converging to $-\infty$ as $\eto$. Other eigenvalues remain bounded as $\eto$. We shall show that these eigenvalues converge to the eigenvalues of $\mathcal{S}_{\alpha,\beta}(\Psi,\Phi)$.

\subsection{Convergence theorem}
Let $\{\lme\}_{\eps\in\mathcal{I}}$ be a sequence of eigenvalues of $\mathcal{S}_\eps$ and assume that $\{\ye\}_{\eps\in\mathcal{I}}$ is a sequence of the corresponding
$L_2(a,b)$-normalized eigenfunctions. Here $\mathcal{I}$ is an infinite subset of $(0,1)$ for which 0 is an accumulation point.

\begin{thm}\label{ConvergenceTheorem}
{\sl
If $\lme\to\lambda$ and $\ye\to v$ in $L_2(a,b)$ weakly as $\mathcal{I}\ni\eto$,
then $\lambda$ is an eigenvalue of $\mathcal{S}_{\alpha,\beta}(\Psi,\Phi)$ with the corresponding eigenfunction $v$.
Furthermore, $\ye$ converges to $v$ in $L_2(a,b)$.
}\end{thm}
We have divided the proof into a sequence of lemmas.
To start with, let us describe the behavior of $\ye$ outside the $\eps$-neighborhood of the origin.
\begin{lem}\label{LemConvergenceC1}
{\sl
Under the assumptions of Theorem~\ref{ConvergenceTheorem},
for every positive $\gamma$ the sequence $\ye$ tends to $v$ weakly in the topology of $W_2^4((a,b)\setminus(-\gamma,\gamma))$ and strongly in the topology of
$C^3([a,b]\setminus(-\gamma,\gamma))$.
Furthermore, $v$ solves the equation
  \begin{equation}\label{UnperturbedSch}
    Lv=\lambda v,\qquad x\in(a,0)\cup(0,b).
   \end{equation}
}\end{lem}
\begin{proof}
Throughout the proof $\mathcal{G}_\gamma$ denotes the set of test functions $\varphi\in C_0^\infty(a,b)$ such that $\varphi(x)=0$ for $x\in (-\gamma,\gamma)$.
From the equation \eqref{MainProblem} for all $\varphi\in \mathcal{G}_\gamma$ and $\eps<\gamma$ we deduce
  \begin{equation}\label{IntegIdentEps}
    \int_a^b L\ye\,\varphi\,dx=\lme\int_a^b \ye\varphi\,dx,
\end{equation}
since $\mathrm{supp}\:\Psi_\eps\subset (-\gamma,\gamma)$.
The right hand side of \eqref{IntegIdentEps} has a limit as $\mathcal{I}\ni\eto$ by assumption, thus
the integral on the left hand side converges for all $\varphi\in \mathcal{G}_\gamma$.
It follows that $\ye \to v$ in $W_2^4((a,b)\setminus(-\gamma,\gamma))$ weakly and thus   \begin{equation*}\label{IntegIdent}
    \int_a^b Lv\,\varphi\,dx=\lambda\int_a^b v\varphi\,dx, \qquad \varphi\in \mathcal{G}_\gamma.
\end{equation*}
From this identity it may be concluded that $v$ solves \eqref{UnperturbedSch} on $(a,b)\setminus(-\gamma,\gamma)$, and so on $(a,0)$ and $(0,b)$, since $\gamma$ is an arbitrary constant.
Applying the imbedding theorem yields convergence of $\ye$ in  $C^3((a,b)\setminus(-\gamma,\gamma))$, which completes the proof.
\end{proof}
We proceed to investigate the behavior of $\ye$ along with its derivatives at the points $x=-\eps$ and $x=\eps$.
\begin{lem}\label{LemmmaYe(eps)}
{\sl
Let $\lambda_\eps\to \lambda$ and $\ye \to v$ in $L_2(a,b)$ weakly as $\mathcal{I}\ni\eto$. Then for $k=0,\ldots,3$ the sequences $\ye^{(k)}(\pm\eps)$ converge to $v^{(k)}(\pm0)$ as
$\mathcal{I}\ni\eps\to 0$.
}
\end{lem}
\begin{proof}
Let $\zeta_k$ be $C^\infty_0((a,b)\setminus0)$-functions such that
$\zeta_k(x)=0$ for $x<0$ and $\zeta_k(x)=\textstyle\frac{x^k}{k!}$ for $x\in (0,\textstyle\frac{b}{2})$.
Denote by $\chi_{(\eps,\infty)}$ the characteristic function of $(\eps,\infty)$ and set $\zeta_k^\eps(x)=\chi_{(\eps,\infty)}(x)\zeta_k(x)$.
We note that $\textstyle\frac{d^{k+1}}{dx^{k+1}}\zeta_k^\eps(x)=0$
for $x\in (\eps,{\textstyle\frac{b}{2}})$.
Multiplying both equalities  \eqref{MainProblem}, \eqref{UnperturbedSch} by $\zeta_0^\eps$ and integrating by parts yield
\begin{align*}
\ye'''(\eps)=-\int^b_{\textstyle\frac{b}{2}}\ye'''\zeta'_0\,dx+\int^b_\eps \ye\zeta_0(V-
\lme)\,dx,\\
v'''(\eps)=-\int^b_{\textstyle\frac{b}{2}}v'''\zeta'_0\,dx+
\int^b_\eps v\zeta_0(V-\lambda)\,dx.
\end{align*}
The right hand sides of the equalities have the same limit as $\mathcal{I}\ni\eto$ in view of Lemma~\ref{LemConvergenceC1}, and so $\ye'''(\eps)\to v'''(+0)$.
Applying the function $\zeta_0^\eps(-x)$ similar to the above implies $\ye'''(-\eps)\to v'''(-0)$. We have proved the Lemma for $k=3$.

The case $k=2$ can be handled in much the same way, the only difference being in applying the function $\zeta_1^\eps$.
Multiplying \eqref{MainProblem}, \eqref{UnperturbedSch} by $\zeta_1^\eps$ and integrating by parts, we derive
\begin{align*}
\ye''(\eps)&=\ye'''(\eps)\eps-\int^b_{\textstyle\frac{b}{2}}\ye''\zeta''_1\,dx+
\int^b_\eps \ye\zeta_1(\lme-V)\,dx,\\
v''(\eps)&=v'''(\eps)\eps-\int^b_{\textstyle\frac{b}{2}}v''\zeta''_1\,dx+
\int^b_\eps v\zeta_1(\lambda-V)\,dx.
\end{align*}
Again employing Lemma~\ref{LemConvergenceC1}, we deduce $u_\eps''(\eps)\to v''(+0)$ as $\mathcal{I}\ni\eto$. Using $\zeta_1^\eps(-x)$ instead of $\zeta_1^\eps(x)$ yields $\ye''(-\eps)\to v''(-0)$.
The rest of the proof runs as before.
\end{proof}

We denote by $g_1$, $g_2$ solutions of the following problems
\begin{align}
\begin{aligned}
 &g_k^{(4)}+\alpha \Psi (\xi)g_k=0,\qquad\xi\in(-1,1),\\\label{CauchyProblemGk}
&g_k(-1)=\delta_{1,k},\qquad g'_k(-1)=\delta_{2,k},\qquad
g''_k(-1)=0,\qquad g'''_k(-1)=0,
\end{aligned}\end{align}
where $\delta_{i,j}$ is the Kronecker symbol.
Let $g_\eps$ solve the Cauchy problem on $[-1,1]$
\begin{align}
\begin{aligned}
 &g^{(4)}+(\alpha \Psi (\xi)+\eps^2\gamma_1\Upsilon_1(\xi))g=-\beta \Phi (\xi)\bigl(\eps^{-1}\ye(-\eps)g_1(\xi)+
 \ye'(-\eps)g_2(\xi)\bigr),\\\label{CauchyProblemGEps}
&g(-1)=0,\qquad g'(-1)=0,\qquad
g''(-1)=v''(-0),\qquad g'''(-1)=0,
\end{aligned}
\end{align}
and let $g_k^\eps$ be the solutions of the problems
\begin{align*}
\begin{aligned}
 &g_k^{(4)}+(\alpha \Psi (\xi)+\eps^2\gamma_1\Upsilon_1(\xi))g_k=0,\qquad\xi\in(-1,1),\\\label{CauchyProblemGk}
&g_k(-1)=\delta_{1,k},\qquad g'_k(-1)=\delta_{2,k},\qquad
g''_k(-1)=0,\qquad g'''_k(-1)=0.
\end{aligned}\end{align*}

The task is now to describe the behavior of the eigenfunction $\ye$ in the $\eps$-neighborhood of the origin.
\begin{lem}\label{LemmaLocalConvergence}
{\sl
If $\lme\to \lambda$ and $\ye \to y$ in $L_2(a,b)$ weakly as $\mathcal{I}\ni\eps\to 0$, then
\begin{equation}\label{EstimYeps0}
    |\ye(-\eps)|+|\ye(\eps)|\leq c\eps,
\end{equation}
and moreover,
\begin{equation}\label{YeG123}
\left\|\eps^{-2}\ye(\eps\xi)-
\eps^{-2}\ye(-\eps)g_1^\eps(\xi)-\eps^{-1}\ye'(-\eps)g_2^\eps(\xi)
- g_\eps(\xi)\right\|_{C^3([-1,1])}\to 0.
\end{equation}
}\end{lem}
\begin{proof}
First let us prove that the sequence $u_\eps(\xi)=\eps^{-1}\ye(\eps\xi)-\eps^{-1}\ye(-\eps)g_1(\xi)-\ye'(-\eps)g_2(\xi)$ tends to 0 in $L_2(-1,1)$ as $\mathcal{I}\ni\eto$.
Set $\mathcal{U}_\eps(\xi)=\eps^{-1}\ye(\eps\xi)-\eps^{-1}\ye(-\eps)h_1^\eps(\xi)-\ye'(-\eps)h_2^\eps(\xi)$ with $h_1^\eps$, $h_2^\eps$ being solutions to problems
\begin{align*}
\begin{aligned}
 &h_k^{(4)}+(\alpha \Psi (\xi)+\eps\beta \Phi (\xi))h_k=0,\qquad\xi\in(-1,1),\\\label{CauchyProblemGk}
&h_k(-1)=\delta_{1,k},\qquad h'_k(-1)=\delta_{2,k},\qquad
h''_k(-1)=0,\qquad h'''_k(-1)=0.
\end{aligned}\end{align*}
By construction $\mathcal{U}_\eps$ solves the problem
\begin{equation*}%\label{CauchyProblemUe}
    \begin{cases}
    u^{(4)}+(\alpha \Psi (\xi)+\eps\beta \Phi (\xi)) u=h_\eps(\xi),\qquad \xi\in(-1,1),\\
    u(-1)=0,\quad     u'(-1)=0,\quad
    u''(-1)=\eps\ye''(-\eps),\quad u'''(-1)=\eps^2 \ye'''(-\eps),
    \end{cases}
\end{equation*}
where $h_\eps(\xi)=-\eps\ye(\eps\xi)(\gamma_1\Upsilon_1(\xi)+\eps\gamma_2\Upsilon_2(\xi)+\eps^2 U(\eps\xi)-\eps^2\lme)$.
For every function $h_\eps\in L_2(-1,1)$ the solution $\mathcal{U}_\eps$ of the above problem is unique, belongs to $W_2^4(-1,1)$, and satisfies the estimate
\begin{equation*}%\label{EstimU}
    \|\mathcal{U}_\eps\|_{W_2^4(-1,1)}\leq c\,\|h_\eps\|_{L_2(-1,1)}
\end{equation*}
with constant $c$ being independent of $\eps$. Employing the inequality
\begin{equation*}
    \eps\int_{-1}^1\ye^2(\eps\xi)\,d\xi=\int_{-\eps}^\eps \ye^2(x)\,dx\leq \|\ye\|_{L_2(a,b)}^2=1,
\end{equation*}
one concludes that $\eps^{1/2}\|\ye(\eps\xi)\|_{L_2(-1,1)}\leq c$, hence that $\|h_\eps\|_{L_2(-1,1)}\leq c\eps^{1/2}$, and finally that $\|\mathcal{U}_\eps\|_{W_2^4(-1,1)}\to0$ as $\mathcal{I}\ni\eto$.
From this it follows that
\begin{align}\label{Conv1}
&    \eps^{-1}\ye(\eps)-\eps^{-1}\ye(-\eps)h^\eps_1(1)-\ye'(-\eps)h^\eps_2(1)\to 0,\\\label{Conv2}
&
\ye'(\eps)-\eps^{-1}\ye(-\eps)(h^\eps_1)'(1)-\ye'(-\eps)(h^\eps_2)'(1)\to 0,\\\label{Conv3}
&    \eps\ye''(\eps)-\eps^{-1}\ye(-\eps)(h^\eps_1)''(1)-\ye'(-\eps)(h^\eps_2)''(1)\to 0,\\\label{Conv4}
&    \eps^2 \ye'''(\eps)-\eps^{-1}\ye(-\eps)(h^\eps_1)'''(1)-\ye'(-\eps)(h^\eps_2)'''(1)\to 0,
\end{align}
By \eqref{Conv2}--\eqref{Conv4} we deduce
\begin{equation*}
\ye(-\eps)(h^\eps_1)^{(k)}(1)=O(\eps)\quad\text{as}\quad\mathcal{I}\ni\eps\to 0\quad\text{for}\quad k=1,2,3.
\end{equation*}
If at least one of the values $(h^\eps_1)''(1)$ or $(h^\eps_1)'''(1)$ is nonzero, then $\ye(-\eps)=O(\eps)$ as $\mathcal{I}\ni\eps\to 0$.
Suppose, to the contrary, that $(h^\eps_1)''(1)=(h^\eps_1)'''(1)=0$. From \eqref{CauchyProblemGk} it follows that $g_1$ is an eigenfunction
of the problem \eqref{w0Problem} corresponding to the eigenvalue $\alpha\in \Sigma_\Psi$.
By construction $g_1'(-1)=0$, contrary to \eqref{MainAssumption}.
The proof of \eqref{EstimYeps0} is complete by using \eqref{Conv1}.

Finally, since $\|g_k-h_k^\eps\|_{L_2(-1,1)}$ converges to $0$ as $\mathcal{I}\ni\eto$ for $k=1,2$, then $u_\eps$ converges to 0 in $L_2(-1,1)$ as $\mathcal{I}\ni\eto$.
In fact,
\begin{multline*}
    \|u_\eps\|_{L_2(-1,1)}\leq  \|\mathcal{U}_\eps\|_{L_2(-1,1)}+\eps^{-1}\ye(-\eps)\|h_1^\eps-g_1\|_{L_2(-1,1)}+\\
    \ye'(-\eps)\|h_2^\eps-g_2\|_{L_2(-1,1)}\to0
\end{multline*}
as $\eto$.

We introduce the function $\mathcal{W}_\eps(\xi)=\eps^{-2}\ye(\eps\xi)-
  \eps^{-2}\ye(-\eps)g_1^\eps(\xi)-\eps^{-1}\ye'(-\eps)g_2^\eps(\xi)
  -g_\eps(\xi)$.
Analysis similar to the above implies that $\mathcal{W}_\eps$ solves the problem
\begin{equation*}%\label{CauchyProblemUe}
    \begin{cases}
    u^{(4)}+(\alpha \Psi (\xi)+\eps^2\gamma_1\Upsilon_1(\xi)) u=f_\eps(\xi),\qquad \xi\in(-1,1),\\
    u(-1)=0,\quad     u'(-1)=0,\quad
    u''=\ye''(-\eps)-v''(-0),\quad u'''(-1)=\eps\ye'''(-\eps),
    \end{cases}
\end{equation*}
with $f_\eps(\xi)=-\beta \Phi u_\eps(\xi)-\eps\ye(\eps\xi)(\gamma_2\Upsilon_2+\eps U(\eps\xi)-\eps\lme)$ and satisfies the estimate
\begin{equation*}%\label{AprioriEst}
    \|\mathcal{W}_\eps\|_{W_2^4(-1,1)}\leq c\left(\|f_\eps\|_{L_2(-1,1)}+|\ye''(-\eps)-v''(-0)|\right)
\end{equation*}
with constant $c$ being independent of $\eps$. According to Lemma~\ref{LemmmaYe(eps)} the right hand side of this estimate converges to zero. Thus $\mathcal{W}_\eps$ tends to 0 in $W_2^4(-1,1)$ as $\mathcal{I}\ni\eto$.
The embedding $W_2^4(-1,1)\subset C^3([-1,1])$ establishes the convergence \eqref{YeG123}.
\end{proof}

Letting $\mathcal{I}\ni\eto$ we conclude that
\begin{align}\label{Convergance1}
&    \eps^{-2}\bigl(\ye(\eps)-\ye(-\eps)g^\eps_1(1)-\eps \ye'(-\eps)g^\eps_2(1)\bigr)-g_\eps(1)\to 0,\\\label{Convergance2}
&    \eps^{-1}\bigl(\ye'(\eps)-\eps^{-1}\ye(-\eps)(g^\eps_1)'(1)-\ye'(-\eps)(g^\eps_2)'(1)\bigr)-g_\eps'(1)\to 0,\\\label{Convergance3}
&    \ye''(\eps)-\eps^{-2}\ye(-\eps)(g^\eps_1)''(1)-\eps^{-1}\ye'(-\eps)(g_2)''(1)-g_\eps''(1)\to 0,\\\label{Convergance4}
&    \eps \ye'''(\eps)-\eps^{-2}\ye(-\eps)(g_1)'''(1)-\eps^{-1}\ye'(-\eps)(g_2)'''(1)-g_\eps'''(1)\to 0,
\end{align}
in light of Lemma~~\ref{LemmaLocalConvergence} for $\xi=1$.

\begin{lem}\label{LemmaL2convergence}
{\sl
Assume that $\lambda_\eps\to \lambda$ and $\ye\to y$ in $L_2(a,b)$ weakly as $\mathcal{I}\ni\eps\to 0$. Then
$\ye\to y$ in $L_2(a,b)$ as $\mathcal{I}\ni\eps\to 0$.
}\end{lem}
\begin{proof}
First we show that $\ye$ is bounded on $[a,b]$ uniformly with respect to $\eps$.
Applying~\eqref{YeG123} and Lemma~\ref{LemmmaYe(eps)}, we see at once that the sequence $\ye$ is uniformly bounded on $[-\eps,\eps]$
\begin{align}\label{Uniform-ee}
    |\ye(\eps\xi)|\leq c\eps^2+
    |\ye(-\eps)g^\eps_1(\xi)|+\eps|\ye'(-\eps)g^\eps_2(\xi)|+
    \eps^2|g_\eps(\xi)|\leq c_1\eps.
\end{align}
Set $\Omega_\eps=(a,b)\setminus (-\eps,\eps)$. Multiplying the equation \eqref{MainProblem} by the function $\chi_{\Omega_\eps}\ye$ and integrating by parts give
\begin{multline*}
   \int\limits_{\Omega_\eps}\ye''^{2}\,dx=
     \int\limits_{\Omega_\eps}(\lambda_\eps-U)\ye^2\,dx-
     \ye'''(-\eps)\ye(-\eps)+\\\ye'''(\eps)\ye(\eps)+
     \ye''(-\eps)\ye'(-\eps)-\ye''(\eps)\ye'(\eps).
\end{multline*}
All terms on the right-hand side are uniformly bounded with respect to $\eps$.
Thus the sequence $\ye$ is bounded in $W^2_2(\Omega_\eps)$, and so in $C^1(\Omega_\eps)$. On account of the above conclusion combining with \eqref{Uniform-ee}, we deduce that $\max\limits_{x\in(a,b)}|\ye(x)|\leq c$ with constant $c$ being independent of $\eps$.

Fix $\gamma>0$. According to Lemma~\ref{LemConvergenceC1} the difference $\ye-y$ has the $L_2(\Omega_\gamma)$-norm less than $\gamma$, provided $\eps$ is small enough. Then
\begin{multline*}
   \|\ye-y\|_{L_2(a,b)}\leq  \|\ye-y\|_{L_2(\Omega_\gamma)}+\|\ye-y\|_{L_2(-\gamma,\gamma)}\leq\\
    \leq\gamma(1+2\max_{x\in(a,b)}|\ye(x)-y(x)|)\leq C\gamma
\end{multline*}
with constant $C$ being independent of $\eps$.
Recall that $\gamma$ may be made arbitrary small, and the proof is complete.
\end{proof}
\medskip

\noindent
\textit{Proof of Theorem~\ref{ConvergenceTheorem}.}
We conclude from Lemmas~\ref{LemConvergenceC1}, \ref{LemmaL2convergence} that $v$ is a solution
of the equation
\begin{equation*}
    Lv=\lambda v,\qquad x\in (a,0)\cup(0,b),
\end{equation*}
satisfies the boundary conditions $v(a)=v'(a)=0$, $v(b)=v'(b)=0$, and  $\norm{v}_{L_2(a,b)}=1$.
Furthermore, $v(0)=0$ according to \eqref{EstimYeps0}, and we are left with the task of showing that $v$ satisfies appropriate coupling conditions at the origin.

Again applying \eqref{EstimYeps0}, we deduce that $\eps^{-1}\ye(-\eps)$ has a limit as
$\mathcal{I}\ni\eto$, which will be denoted by $s$, and $\ye'(-\eps)\to v'(-0)$ as
shown in Lemma~\ref{LemmmaYe(eps)}.
Therefore $q_\eps=\eps^{-1}\ye(-\eps)g^\eps_1+\ye'(-\eps)g^\eps_2$ converges to $q=s g_1+v'(-0)g_2$ in $C^3([-1,1])$ as $\mathcal{I}\ni\eto$.
From \eqref{Convergance3}--\eqref{Convergance4} it may be concluded that
sequences $\eps^{-1}q_\eps''(1)$ and $\eps^{-1}q_\eps'''(1)$ are bounded as
$\mathcal{I}\ni\eto$, hence that
\begin{equation}\label{Co1}
    q''(1)=0,\qquad q'''(1)=0.
\end{equation}
Combining \eqref{CauchyProblemGk} with \eqref{Convergance2} yields
\begin{equation}\label{Co2}
    q'(-1)=v'(-0),\qquad q'(1)=v'(+0).
\end{equation}
We see at once that $q$ is a solution of the Cauchy problem
\begin{equation}\label{ProblemQ}
    \begin{cases}
    q^{(4)}+\alpha \Psi  q=0,\qquad \xi\in(-1,1),\\
    q(-1)=s,\quad q'(-1)=v'(-0),\quad
    q''(-1)=0,\quad q'''(-1)=0.
    \end{cases}
\end{equation}
Coupling conditions of the limit problem depend on whether the problem \eqref{ProblemQ} admits a nontrivial solution.
Let us suppose for the moment that the problem \eqref{ProblemQ} has a trivial solution $q=0$ only. Then \eqref{Co2} implies the coupling condition $v'(0)=0$.

Next assume that \eqref{ProblemQ} has a nontrivial solution. By \eqref{Co1} $q$ is an eigenfunction of the problem~\eqref{w0Problem} and $\alpha $ belongs to the resonant set $\Sigma_{\Psi }$.
In view of \eqref{Co2} we have
\begin{equation*}
    q'(1)v'(-0)-q'(-1)v'(+0)=0,
\end{equation*}
which is equivalent to
$v'(+0)-\theta_{\Psi }(\alpha )v'(-0)=0$.

 For every
$\eps>0$ the Lagrange identity holds
\begin{equation*}%\label{ExistanceConditionQ1}
   q_\eps'''(1)q(1)-q_\eps''(1)q'(1)=0.
\end{equation*}
Dividing the above identity by $\eps$ and letting $\mathcal{I}\ni\eto$, we derive
\begin{equation}\label{ExistanceCondition}
   z'''(1)q(1)+(v''(+0)-z''(1))q'(1)=0.
\end{equation}
in light of \eqref{Convergance3}, \eqref{Convergance4}, where $z$ solves the problem
\begin{equation*}
\begin{cases}
z^{(4)}+\alpha \Psi (\xi)z=-\beta \Phi (\xi)q(\xi),\qquad \xi\in(-1,1),\\
z(-1)=0,\quad z'(-1)=0,\quad
z''(-1)=v''(-0),\quad z'''(-1)=0.
\end{cases}\end{equation*}
Taking into account \eqref{ExistanceCondition}, the Lagrange identity for $z$ may be written as
\begin{multline*}
z'''(1)q(1)-z''(1)q'(1)+v''(-0)q'(-1)=
\big( z'''(1)q(1)+(v''(+0)-z''(1))q'(1)\big)-\\-v''(+0)q'(1)+v''(-0)q'(-1)=
-v''(+0)q'(1)+v''(-0)q'(-1)
=-\beta \vartheta_\Psi[q].
\end{multline*}
Dividing the last equality by $q'(-1)$ and recalling \eqref{Co2} gives
$$    \theta_{\Psi }(\alpha ) v''(+0)-v''(-0)-\beta \vartheta_\Psi\big[w_{\alpha }/w'_{\alpha }(-1)\big] v'(-0)=0.
$$
Thus $v$ is an eigenfunction of $\mathcal{S}_{\alpha,\beta}(\Psi,\Phi)$ corresponding to the eigenvalue $\lambda$.
\hfill $\square$

\medskip
Theorem~\ref{ConvergenceTheorem} allows one to justify the choice of $\mathcal{S}_{\alpha,\beta}(\Psi,\Phi)$.
\begin{thm}\label{MainTheorem}
{\sl
Suppose that the eigenvalue $\lme$ of $\mathcal{S}_\eps$ is bounded
from below. Then $\lme$ has a finite limit as $\mathcal{I}\ni\eto$ and this limit is
a point of the spectrum of $\mathcal{S}_{\alpha,\beta}(\Psi,\Phi)$.
For each simple eigenvalue  $\lambda$ of $\mathcal{S}_{\alpha,\beta}(\Psi,\Phi)$ there
exist exactly one eigenvalue $\lme$ of $\mathcal{S}_\eps$ converging to $\lambda$ as $\eto$. }
\end{thm}
\begin{proof}
Suppose to start with that
$$
\mu_*=\underline{\lim}_{\eto}\lme_k<
\overline{\lim}_{\eto}\lme_k=\mu^*.
$$
The constants $\mu_*$, $\mu^*$  are finite since $\lme_k$ is a bounded function. Recall that $\lme_k$ is a continuous function of $\eps\in(0,1)$. Then for each $\lambda\in[\mu_*,\mu^*]$ there exists a subsequence of eigenvalues $\lme$, $\eps\in\mathcal{I}$, converging to $\lambda$.The sequence $\{\ye\}_{\eps\in\mathcal{I}}$ of the corresponding normalized eigenfunctions contains a weakly convergent subsequence. By Theorem~\ref{ConvergenceTheorem}, $\lambda$ is an eigenvalue of $\mathcal{S}_{\alpha,\beta}(\Psi,\Phi)$.
Therefore the interval $[\mu_*,\mu^*]$ belongs to the spectrum $\sigma(\mathcal{S}_{\alpha,\beta}(\Psi,\Phi))$, a contradiction.

We now turn to the second part of the theorem.
Let us assume that $\lme_k\to\lambda$ and $\lme_{k+1}\to\lambda$ for some $k$.
Then there exist two sequences $\{y_k^\eps\}_{\eps\in\mathcal{I}}$ and $\{y_{k+1}^\eps\}_{\eps\in\mathcal{I}}$ of eigenfunctions,
which converge in $L_2(a,b)$ to vectors of the form $e^{i\varphi}v$. This contradicts
the fact that $y_k^\eps$ and $y_{k+1}^\eps$ are orthogonal in $L_2(a,b)$ for all $\eps\in\mathcal{I}$.
\end{proof}
\subsection{Approximation theorem}
We proceed to show that each  point of $\sigma(\mathcal{S}_{\alpha,\beta}(\Psi,\Phi))$ is a limit of
the eigenvalues of  $\mathcal{S}_{\varepsilon}$.

Let $B$ be a self-adjoint operator in a Hilbert space  $H$ with domain
$\mathcal{D}(B)$. A pair $(\mu, u)\in \mathbb{R}\times\mathcal{D}(B)$ with $\|u\|_H=1$
is called a \emph{quasimode} of the operator $B$ with an accuracy $\rho>0$ if $\|Bu-\mu u\|_H\leq \rho$.

\begin{lem}\label{LemmaVishik}
{\sl
Suppose that the spectrum of $B$ is discrete and simple. If $(\mu, u)$ is a quasimode of $B$ with accuracy $\rho>0$, then the interval  $[\mu-\rho,\mu+\rho]$ contains an eigenvalue $\lambda$ of $B$.
Furthermore, if the segment $[\mu-\tau,\mu+\tau]$ contains only this eigenvalue of $B$, then $\|u-v\|_H\leq 2\tau^{-1}\rho$, where $v$ is a normalized eigenfunction of $B$ for the eigenvalue $\lambda$.
\rm{\cite{VishykLust}}
}\end{lem}

Let us construct the quasimodes of $\mathcal{S}_{\eps}$.
Suppose $\lambda$ is a simple eigenvalue of the operator $\mathcal{S}_{\alpha,\beta}(\Psi,\Phi)$ with the eigenfunction $v$ such that $\|v\|=1$. Here and subsequently, $\norm{\cdot}$ stands for the norm in $L_2(a,b)$.  For each $\lambda$ and $v$ we have obtained the formal asymptotic approximations
$\Lambda_\varepsilon$, $Y_\varepsilon$ defined by either \eqref{Approx} or \eqref{Approx1} depending on  $\alpha $ and $\Psi $.
In further computation we do not distinguish the resonant and non-resonant cases. By construction we have
\begin{align}
\begin{aligned}\label{Reminders}
    &LY_\eps-\Lambda_\eps Y_\eps=\eps^3 R_1(\eps,x),\qquad x\in(a,-\eps)\cup(\eps,b),\\
    \big(&L+\alpha \eps^{-4}\Psi (\eps^{-1}x)+\beta \eps^{-3}\Phi (\eps^{-1}x)+
    \gamma_1\eps^{-2}\Upsilon_1(\eps^{-1}x)+\gamma_2\eps^{-1}\Upsilon_2(\eps^{-1}x)
    \big) Y_\eps-\\
    &\qquad\qquad\qquad\qquad\qquad\qquad\qquad\qquad\qquad\qquad -\Lambda_\eps Y_\eps=\eps R_2(\eps,x),\qquad x\in(-\eps,\eps).
    \end{aligned}
\end{align}
The function $Y_\eps$ does not belong to the domain of $\mathcal{S}_\eps$, since it has jump dis\-con\-ti\-nuities at the points $\pm\eps$. Indeed,
\begin{multline*}
    [Y_\eps]_{x=\pm\eps}=\eps^3 r_0^\pm(\eps),\quad
    [Y'_\eps]_{x=\pm\eps}=\eps^3 r_1^\pm(\eps),\\
    [Y''_\eps]_{x=\pm\eps}=\eps^3 r_2^\pm(\eps),\quad
    [Y'''_\eps]_{x=\pm\eps}=\eps^2 r_3^\pm(\eps).
\end{multline*}
Here all the functions $R_j$, $r_j^\pm$ are uniformly bounded with respect to their arguments.
We can construct a function $\zeta_\varepsilon$ with the following properties
\begin{itemize}
  \item [$\circ$]$\zeta_\varepsilon$ is a smooth function outside the points $x=\pm\varepsilon$ and differs from zero only for $\varepsilon<|x|<1$;
  \item [$\circ$]$[\zeta_\eps]_{x=\pm\eps}=-\eps r_0^\pm(\eps),$ $[\zeta'_\eps]_{x=\pm\eps}=-\eps r_1^\pm(\eps),$
      $[\zeta''_\eps]_{x=\pm\eps}=-\eps r_2^\pm(\eps)$ and
            $[\zeta'''_\eps]_{x=\pm\eps}=-r_3^\pm(\eps);$
  \item [$\circ$] $\max\limits_{\varepsilon<|x|<1}(|\zeta_\eps(x)|+|\zeta'_\eps(x)|+
      |\zeta''_\eps(x)|+|\zeta'''_\eps(x)|+|\zeta''''_\eps(x)|)\leq c$ with constant $c$ being inde\-pen\-dent of $\varepsilon$,
\end{itemize}
which eliminates discontinuity. In fact, $Y_\varepsilon+\varepsilon^2\zeta_\varepsilon$ is a function from $C^3(a,b)$ and belongs to $\mathcal{D}(\mathcal{S}_{\varepsilon})$. We set $\mathcal{Y}_\varepsilon=\|Y_\varepsilon+
\varepsilon^2\zeta_\varepsilon\|^{-1}(Y_\varepsilon+\varepsilon^2\zeta_\varepsilon)$
and substitute $\mathcal{Y}_\varepsilon$ into \eqref{Reminders} instead of $Y_\varepsilon$. Then
the orders of smallness of  right-hand sides in \eqref{Reminders} do not  change since $\|Y_\varepsilon+
\varepsilon^2\zeta_\varepsilon\|\to 1$ as $\varepsilon\to 0$.
Therefore the  pair $(\Lambda_\varepsilon, \mathcal{Y}_\varepsilon)$ is a quasimode of  $\mathcal{S}_{\varepsilon}$ with accuracy $\varepsilon$.
\begin{thm}
{\sl
Given $(\alpha,\beta,\gamma_1,\gamma_2;\Psi,\Phi,\Upsilon_1,\Upsilon_2)
\in\mathbb{R}^4\times (C_0^\infty(-1,1))^4$, suppose that $\lambda$ is a simple eigenvalue of $\mathcal{S}_{\alpha,\beta}(\Psi,\Phi)$ with the normalized eigenfunction $v$. Then there exists a simple eigenvalue $\lambda^\eps_j$ of $\mathcal{S}_\eps$ with a corresponding normalized eigenfunction $y^\eps_j$ such that
\begin{align}\label{EstLambda}
\abs{\lambda^\eps_j-\lambda}\leq c_1\eps,\qquad
\bigl\|y^\eps_j-v\bigr\|\leq c_2\eps
\end{align}
with constants $c_1$, $c_2$ being independent of $\eps$.
}\end{thm}

\begin{proof}
Let $(\Lambda_\eps,\mathcal{Y}_\eps)$ be a quasimode of $\mathcal{S}_\eps$ corresponding to the limit eigenvalue $\lambda$ and the eigenfunction $v$. According to Lemma~\ref{LemmaVishik} there exists an eigenvalue $\lambda^\eps_j$ such that $|\lambda^\eps_j-\Lambda_\eps|\leq c_1\eps$, from which the first inequality in \eqref{EstLambda} follows.
In view of Theorem~\ref{MainTheorem} the index $j$ is independent of $\eps$.
If $\tau$ is less than the distance from $\lambda$ to the rest of the spectrum of $\mathcal{S}_{\alpha,\beta}(\Psi,\Phi)$, then the interval $[\lambda-\tau,\lambda+\tau]$ contains the eigenvalue $\lme_j$ only, provided $\eps$ is small enough. Applying again Lemma~\ref{LemmaVishik} yields $\norm{y^\eps_j-\mathcal{Y}_\eps}\leq2\tau^{-1}c_1\eps$, from which the second inequality in \eqref{EstLambda} immediately follows.
\end{proof}

% ----------------------------------------------------------------


\begin{thebibliography}{99}
\bibitem{Alb2ed}
 S. Albeverio,  F. Gesztesy,  R. H{\o}egh-Krohn,  H. Holden,
Solvable models in quantum mechanics. With an appendix by Pavel Exner,
RI: AMS Chelsea Publishing, 2005.

\bibitem{AlbKur}
 S. Albeverio, P. Kurasov,  Singular perturbations of differential operators
and solvable Schr\"{o}dinger type operators,  Cambridge, Univ. Press, 2000.

\bibitem{IA}
T. Azizov, I. Iokhvidov,  Linear operators in space with an indefinite metric, Pure and Applied Mathematics, Chichester, 1989.

\bibitem{BirmanSolom}
M.~Birman, M.~Solomyak, {\it Estimates for the number of negative eigenvalues of the Schr\"odinger operator and its generalization,}  Adv. Sov. Math. {\bf 7} (1991), 1--55.

\bibitem{BerShu} F. Berezin, M. Shubin, The Schr\"{o}dinger equation, Kluwer Academic
Publishers, 1991.

\bibitem{BerezinFadeev}
F. Berezin, L. Faddeev, {\it A remark on Schr\"odinger's equation with a singular potential,} Sov. Math. Dokl. {\bf 2} (1961), 372--375.

\bibitem{CurJDE}
 B. \'{C}urgus, H. Langer,
{\it A Krein space approach to symmetric ordinary differential operators with
    an indefinite weigth function,} J. Diff. Eq. {\bf 79} (1989), 31--61.

\bibitem{DemkovOstrovskii}
Yu.~Demkov, V.~Ostrovskii, Zero-range potentials and their applications in
atomic physics, Leningrad, Leningrad Univ. Press, 1975.

\bibitem{GolshOri}
N. Goloscshapova, L. Oridoroga, {\it 4-th order differential operator with local point
interactions,} Ukr. Math. Bull. {\bf4} (2007), 355--369.

\bibitem{GolovManko0}
Yu. Golovaty, S. Man'ko, {\it Schr\"{o}dinger operator with $\delta'$-potential,} Dopov. Nats. Akad. Nauk Ukr., Mat. Pryr. Tekh. Nauky. {\bf5} (2009), 16--21.

\bibitem{GolovManko1}
Yu. Golovaty, S. Man'ko, {\it Solvable models for the Schr\"{o}dinger operators with $\delta'$-like potentials,} Ukr. Mat. Visn. {\bf 6} (2009), 173--206.

\bibitem{Manko}
S. Man'ko,\textit{ Fourth order differential operators with distributions in coefficients,} Visn. L'viv. Univ., Ser. Mekh.-Mat. {\bf 72} (2010)
(accepted for publication).

\bibitem{Na}
M. Najmark, Linear differential operators, Moscow, Nauka, 1969.







%\bibitem{GoloschapovaOridoroga}
%N. I. Goloscshapova, L. L. Oridoroga.\textit{ 4-th order differential operator
%with local point interactions }// Ukrainian Mathematical Bulletin. 4 (2007),
%N 3, 355--369.


 \bibitem{VishykLust}
M. Vishik, L. Lyusternik,
{\it  Regular degeneration and boundary layer for linear differential
    equations with small parameter,} Usp. Mat. Nauk {\bf12} (1957), 3--122.



\bibitem{Yablonski}
A. Yablonski, {\it Differential equations with generalized coefficients,}
Nonlinear Anal., Theory Methods Appl. {\bf63} (2005), 171--197.

\bibitem{Zavalishchin}
S.~Zavalishchin, A.~Sesekin, Dynamic impulse systems: theory and applications,
Dordrecht, Kluwer Academic Publishers Group, 1997.

\end{thebibliography}
\end{document}